\documentclass[leqno,12pt]{article}
\usepackage{amsmath,amssymb,a4wide,color}

\usepackage[curve,matrix,arrow,cmtip]{xy}

\usepackage{verbatim}

\NoComputerModernTips

\newdir^{ (}{{}*!/-3pt/\dir^{(}}
\newdir_{ (}{{}*!/-3pt/\dir_{(}}

\def\OM{\mathchoice
  {\rlap{\kern3.2pt$\overline{\phantom{L}}$}M}
  {\rlap{\kern3.2pt$\overline{\phantom{L}}$}M}
  {\rlap{\kern2.4pt$\scriptstyle\overline{\phantom{L}}$}M}
  {\rlap{\kern1.8pt$\scriptscriptstyle\overline{\phantom{L}}$}M}}

\let\le\leqslant
\let\ge\geqslant

 \def\End{\mathop{\rm End}\nolimits}

 \def\deg{\mathop{\rm deg}\nolimits}

\def\GL{\mathop{\rm GL}\nolimits}

\let\phi\varphi
\let\theta\vartheta
\let\epsilon\varepsilon

\newtheorem{Thm}{Theorem}

\newtheorem{Lem}[Thm]{Lemma}

\newtheorem{Cor}[Thm]{Corollary}

\newtheorem{Def}[Thm]{Definition}

\newtheorem{Rem}[Thm]{Remark}

\def\qed{{\hskip0pt\unskip\unskip\nobreak\hfil\penalty50
          \hskip1em\hbox{}\nobreak\hfil
           {$\square$}
          \parfillskip=0pt\finalhyphendemerits=0
          \par}\medskip}

\newenvironment{Proof}
               {\noindent{\bf Proof.}\ }
               {\qed}

\newenvironment{Proofof}[1]
               {\noindent{\bf Proof of #1.}\ }
               {\qed}

\newcommand{\BC}{{\mathbb{C}}}

\newcommand{\BF}{{\mathbb{F}}}

\newcommand{\BN}{{\mathbb{N}}}

\newcommand{\BP}{{\mathbb{P}}}

\newcommand{\BZ}{{\mathbb{Z}}}

\newbox\mybox
\def\arrover#1{\mathrel{
       \setbox\mybox=\hbox spread 1.4em
              {\hfil$\scriptstyle#1$\hfil}
       \vbox{\offinterlineskip\copy\mybox
             \hbox to\wd\mybox{\rightarrowfill}}}}

\def\larrover#1{\mathrel{
       \setbox\mybox=\hbox spread 1.4em
              {\hfil$\scriptstyle#1\vphantom{g}$\hfil}
       \vbox{\offinterlineskip\copy\mybox
             \hbox to\wd\mybox{\leftarrowfill}}}}

\def\ontoover#1{\mathrel{
       \setbox\mybox=\hbox spread 1.4em
              {\hfil$\scriptstyle#1\vphantom{g}$\hfil}
       \vbox{\offinterlineskip\copy\mybox
             \hbox to\wd\mybox{\rightarrowfill\hskip-2.8mm
                               $\rightarrow$}}}}
\def\leftontoover#1{\mathrel{
       \setbox\mybox=\hbox spread 1.4em
              {\hfil$\scriptstyle#1\vphantom{g}$\hfil}
       \vbox{\offinterlineskip\copy\mybox
             \hbox to\wd\mybox{$\leftarrow$\hskip-2.8mm
                               \leftarrowfill}}}}

\def\Cinf{{\BC}_\infty}
\def\Finf{F_\infty}

\newbox\dotDdbox
\newbox\dotDtbox
\newbox\dotDsbox
\newbox\dotDssbox
\setbox\dotDdbox=\vbox{\hbox{\kern3.5pt\bf.}\vskip-4.5pt\hbox{$D$}}
\setbox\dotDtbox=\vbox{\hbox{\kern3.5pt\bf.}\vskip-4.5pt\hbox{$D$}}
\setbox\dotDsbox=\vbox{\hbox{\kern2.1pt.}
                       \vskip-6.6pt\hbox{$\scriptstyle D$}}
\setbox\dotDssbox=\vbox{\hbox{\kern1.6pt.}
                        \vskip-8.1pt\hbox{$\scriptscriptstyle D$}}

\newcommand{\prodprm}{\sideset{}{'}\prod}

\newcommand{\abs}[1]{\left\lvert #1\right\rvert}
\newcommand{\period}{\xi}

\begin{document}

\title{A product formula for the higher rank Drinfeld discriminant function}

\author{Dirk Basson\footnote{This work
was supported by the Wilhelm Frank Scholarship.}\\
{\small Department of Mathematical Sciences, 
Stellenbosch University,} \\
{\small Private Bag X1,
Matieland, 7602,
South Africa}\\
{\small djbasson@sun.ac.za}}

\date{}
\maketitle


\begin{abstract}
We give a product expansion for the Drinfeld discriminant function in arbitrary rank $r$, 
which generalizes the formula obtained by Gekeler for the rank 2 Drinfeld discriminant 
function. This enables one to compute the Fourier expansion  of this function much 
more efficiently.

The formula in this article uses an $r-1$-dimensional parameter 
and as such provides a nice counterpoint to the formula previously obtained by Hamahata, 
which is written in terms of several 1-dimensional parameters.
\end{abstract}


\section{Introduction}
One of the first great theorems in the theory of elliptic modular forms is the 
Jacobi product formula for the discriminant function:
\[ \Delta(\tau) = (2\pi i)^{12} q\prod_{n\ge 1}(1-q^n)^{24},\]
where $q=e^{2\pi i\tau}$.  After Drinfeld introduced elliptic modules \cite{Drinfeld1}, 
the study of a characteristic $p$ version of modular forms, called
Drinfeld modular forms, took off beginning with David Goss's  
thesis \cite{GossES}. It was not long before Gekeler \cite{GekelerProd} proved a version of Jacobi's product 
formula for Drinfeld modular forms of rank 2.
 But the higher rank case had to 
wait for some new ideas to make a higher dimensional theory of modular forms 
possible. This came in the form of Pink's new compactification \cite{Pink} of Drinfeld moduli 
spaces with so-called fine level structure which built on a construction of Kapranov in a 
special case \cite{Kapranov}. A project initiated by Breuer and Pink 
and carried on by the author in his PhD dissertation translated this idea into 
the rigid analytic setting \cite{BBP}. It is fitting that a generalization of the product 
formula (see Theorem \ref{Thm:prodformula}) is one of the first results proved in the theory of higher rank Drinfeld 
modular forms.

It should be noted that there exists a product formula for higher rank Drinfeld 
discriminant functions, discovered by Hamahata (\cite{Hamahata} and \cite{Hamahata2}). 
It has a very different flavour from the formula proven in this note, since it uses 
$r$ parameters each depending on one variable, while our formula uses one parameter
that depends on $r-1$ variables.

Generally speaking, expansions of Drinfeld modular forms are quite hard to compute. This 
is true for rank 2 Drinfeld modular forms, as can be seen from \cite{GekelerCoeff}, where
even the simplest functions have complicated expansions; 
and things only become worse in higher rank when 
coefficients are no longer constants, but functions. The 
product formula
simplifies the computation of expansions for higher rank Drinfeld modular forms
and is therefore not just an interesting result, but a useful one as well.\\

\section{Background and notation}
Throughout this article, $A=\BF_q[t]$ and $F=\BF_q(t)$. Also let $\Finf=\BF_q((1/t))$ be the 
completion at the infinite place for which $\abs{f/g}=q^{\deg f-\deg g}$ and let $\Cinf$ denote 
the completion of an algebraic closure of $\Finf$.
The Drinfeld period domain $\Omega^r$ is the complement of all the $\Finf$ hyperplanes in 
$\BP^{r-1}(\Cinf)$. It inherits the structure of a rigid analytic space from 
$\BP^{r-1}(\Cinf)$ since it is an admissible open subset \cite[Proposition 1]{SchneiderStuhler}. 
We shall need a more explicit
description of this structure before the technical Lemma \ref{lem:UniformConvergence}, 
but this will be postponed until it is needed.

We shall write the elements in $\Omega^r$ as $r\times 1$ column vectors and normalize them 
so that the last entry of $\omega=[\omega_1,\omega_2,\ldots,\omega_{r-1},\omega_r]^T$ is 
equal to $\period$, where 
\[\period=\sqrt[q-1]{-t^q}\cdot \prod_{i\ge 1}\left(1-\frac{t^{q^i}-t}{t^{q^{i+1}}-t}\right)\]
 is a generator for the 
lattice corresponding to the Carlitz module. (This is only defined up to multiplication by 
an element in $\BF_q^*$, but we choose one and fix it.)
By definition of $\Omega^r$, the elements $\omega_1,\ldots,\omega_r$
are $\Finf$-linearly independent and, in particular, non-zero. 

There is an action of $\GL_r(F)$ on $\Omega^r$ given by 
\[ \gamma\cdot\omega = j(\gamma,\omega)^{-1}[\gamma][\omega], \]
where $[\gamma][\omega]$ denotes the matrix product of the $r\times r$ matrix $[\gamma]$ with the $r\times 1$ matrix $[\omega]$;
 and $j(\gamma,\omega)$ is a 
factor used to normalize $\gamma\cdot\omega$ so that its last entry is $\period$. 
Explicitly, $j(\gamma,\omega)$ is $\period^{-1}$ times the last entry of 
$[\gamma][\omega]$. \\

If an $\BF_q$ subvector space $L\subset \Cinf$ has the property that its intersection
with any ball of finite radius is finite, then we can associate an exponential function
$e_L:\Cinf\to\Cinf$ to it. Under the convention that a primed sum or product indicates
a sum or product over all non-zero elements of a set, it is defined as
\[e_L(X) = X\prodprm_{\lambda\in L}\left(1-\frac{X}{\lambda}\right),\]
and is an entire surjective function with $L$ as its set of roots. It is also $\BF_q$-linear
in the sense that
\[ e_L(ax+by) = ae_L(x) + be_L(y)\]
for all $x,y\in \Cinf$ and $a,b\in \BF_q$. It is easily verified that simultaneously 
scaling the lattice and the parameter $x$ by a factor $c\in \Cinf$ has the following effect:
\[ e_{cL}(cx) = c e_L(x).\]

If $L$ is moreover an $A$-module, then there is the additional identity
\[ e_L(ax) = \varphi^L_a\left(e_L(x)\right), \]
where $\varphi^L_a$ is the polynomial
\begin{equation} \label{eq:DrinfModDef}
\varphi^L_a(Z) = aZ\prodprm_{\lambda\in a^{-1}L/L} \left(1-\frac{Z}{e_L(\lambda)}\right).
\end{equation}
(Note that since $e_L$ is $L$-invariant, the value of $e_L(\lambda)$ does not depend on which representative in $a^{-1}L$ 
is chosen.)
Each $\varphi^L_a$ can be thought of as an $\BF_q$-linear endomorphism of $\Cinf$, and then
the map $\varphi^L:A\to\End_{\BF_q}(\Cinf)$ defined by $a\mapsto\varphi^L_a$ actually defines a 
ring homomorphism. We call this homomorphism the Drinfeld module associated to $L$. 
For proofs of these results and more about Drinfeld modules, see \cite{GossBS}.\\

Any $\omega\in\Omega^r$ gives rise to a lattice $\Lambda_\omega:= A^r\omega\subset \Cinf$
and hence also to a Drinfeld module which we denote $\varphi^\omega$ instead of $\varphi^{ A^r\omega}$.
%
The leading coefficient of the polynomial $\varphi^\omega_t(X)$ depends 
on $\omega$, and in fact defines a rigid analytic function $\Delta:\Omega^r\to \Cinf$, 
which is the function of interest in this note.

If we act with $\gamma\in\GL_r(A)$ on $\omega$, then the lattice $A^r\omega $ maps to 
$\Lambda_{\gamma\omega}=A^r(\gamma\cdot\omega)=j(\gamma,\omega)^{-1}A^r\omega $. In particular
\begin{align*}
 \varphi^{\gamma\cdot\omega}_t(X) 
&= tX\prodprm_{\lambda\in t^{-1}\Lambda_{\gamma\omega}/\Lambda_{\gamma\omega}}\left(1-\frac{X}{e_{\Lambda_{\gamma\omega}}(\lambda)}\right)\\
&\stackrel{\eqref{eq:DrinfModDef}}{=\joinrel=} tX\prodprm_{\lambda\in t^{-1}\Lambda_\omega/\Lambda_\omega}
  \left(1-\frac{X}{j(\gamma,\omega)^{-1}e_{\Lambda_\omega}(\lambda)}\right).
\end{align*}
This equation shows that the leading coefficient changes according to the formula
\[ \Delta(\gamma\cdot \omega)=j(\gamma,\omega)^{q^r-1}\Delta(\omega).\]
In other words $\Delta$ is a (weak) modular form of weight $q^r-1$:

\begin{Def}
Let $k$ be an integer. 
A \emph{weak Drinfeld modular form of rank $r$ and weight $k$} is a rigid analytic function 
$f:\Omega^r\to\Cinf$ satisfying $f(\gamma\cdot\omega)=j(\gamma,\omega)^k f(\omega)$ 
for all $\omega\in\Omega^r$ and all $\gamma\in\GL_r(A)$.
\end{Def}

As in the classical and rank 2 cases, weak modular forms also admit a Fourier expansion. 
Write $\omega'=[\omega_2,\omega_3,\ldots,\omega_r]^T$ for the element of $\Omega^{r-1}$
obtained from $\omega$ by deleting its first component. 
We introduce a parameter $u_{\omega'}(\omega_1)$ defined as
\[ u_{\omega'}(\omega_1) = e_{ A^{r-1}\omega'}(\omega_1)^{-1}.\] 
We shall often abbreviate $u:= u_{\omega'}(\omega_1)$.

\begin{Thm}[\cite{BBP}]
Every weak modular form $f$ has a Laurent series expansion
\[ f(\omega) = \sum_{n\in\BZ} f_n(\omega')u_{\omega'}(\omega_1)^n, \]
where each $f_n:\Omega^{r-1}\to\Cinf$ is a holomorphic function. Moreover, the $f_n$
are uniquely determined.
\end{Thm}

\begin{Rem}
The reader may be interested to know that the functions $f_n$ are themselves also
weak Drinfeld modular forms of rank $r-1$ and weight $k-n$. This fact is proven in \cite{BBP}.
\end{Rem}

\begin{Def}
Let $k$ be an integer. 
A \emph{modular form $f$  of rank $r$ and weight $k$} is a weak modular form  of rank 
$r$ and weight $k$ whose $u$-expansion has the property that $f_n$ is the zero function 
for all $n<0$.
\end{Def}
That concludes the background we need from \cite{BBP}.

\section{Main result and proof}
From now on we denote the lattices $A^r\omega $ and $A^{r-1}\omega' $ by $\Lambda$ and $\Lambda'$
respectively and assume that the Drinfeld module $\varphi$ is the one 
associated to $\Lambda'$. Suppose that the leading coefficient of
the polynomial $\varphi_a(X)$ is $\Delta'_a(\omega')$, where the argument stresses the fact
that this is a function depending on $\omega'$.
For each $a\in A$ we define
\[
f_a(X)=X^{q^{(r-1)\deg a}}\Delta'_a(\omega')^{-1}\varphi_a(X^{-1})-1.
\]
It is a polynomial of degree $q^{(r-1)\deg a}-1$ which is divisible by 
$X^{q^{(r-1)\deg a}-q^{(r-1)\deg a-1}}$.
Then
\begin{equation}
\varphi_a(u^{-1})^{-1}=u^{q^{(r-1)\deg a}}\Delta'_a(\omega')^{-1}(1+f_a(u))^{-1}.
\label{eq:UaInTermsOfHa}
\end{equation}
Note for later use that
\begin{equation}\label{eq:DeltaA}
\Delta'_a(\omega')=\Delta'_t(\omega')^{1+q+\cdots+q^{\deg a-1}},
\end{equation}
which can be obtained by computing $\varphi_a(X)$ from $\varphi_t(X)$. In the proof of
Theorem \ref{Thm:prodformula} we suppress the 
subscript when we mean $\Delta'=\Delta'_t$.

The product in Lemma \ref{lem:UniformConvergence} below is closely related to the product in the expansion 
of $\Delta$. We need to deal with its convergence issues as some point, so we start with it. 
Then we may do the various interchanges of summations and/or products in the proof without worry.
 Lemma \ref{lem:UniformConvergence} asserts that the product converges uniformly on a certain set, which we need to describe
first.

Since $\Omega^r$ is the complement of $\Finf$ hyperplanes in $\BP^{r-1}(\Cinf)$, we wish to 
have a measure of how far an element $\omega\in\Omega^r$ is from any hyperplane. For this we follow the approach
in \cite{SchneiderStuhler}. If $H$ is a 
hyperplane defined over $\Finf$, it can be defined as the zero set of a linear equation (with coefficients in $\Finf$)
\[ \ell_H(x) = h_1x_1+h_2x_2+\cdots+h_rx_r = 0.\]
This linear equation is only well-defined up to scaling by a factor in $\Cinf^\times$, so we
assume that $\max{\abs{h_j}}=1$. Then, even though $\ell_H(\omega)$ might not be well-defined, $\abs{\ell_H(\omega)}$ is well-defined for any hyperplane $H$ and any $\omega\in \Omega^r$. We define 
\[ \abs{\omega}_i = \inf\{\abs{\ell_H(\omega)} : H\text{ is a hyperplane defined over }\Finf\}. \]
If we also define $\abs{\omega}= \max_{j=1}^r \abs{\omega_j}$, then the sets
\[\Omega^r_n := \{\omega\in\Omega^r\mid \abs{\omega}_i \ge q^{-n}\abs{\omega}\}\]
for $n\in\BN$ form a collection of affinoid subdomains of $\Omega^r$.

\begin{Rem}
Let $B(0,r)=\{z\in\Cinf\mid \abs{z}<r\}$ be the open ball with radius $r$ in $\Cinf$. The sets 
\[ \bigcup_{n\ge 1} B(0,r_n)\times \Omega^{r-1}_n \]
(where intuitively $r_n$ gets smaller as $n$ increases) play an important role in the 
analysis of higher rank Drinfeld modular forms and are called
\emph{quasi-uniform neighbourhoods} (of $\{0\}\times\Omega^{r-1}$) in \cite{BBP}. 
Incidentally, the set on which the $u$-expansion of a Drinfeld modular form converges
is a quasi-uniform neighbourhood (punctured to exclude $\{0\}\times\Omega^{r-1}$ if 
it is not holomorphic at infinity). This at least gives some motivation as to why it 
occurs in the next Lemma.
\end{Rem}

\begin{Lem}\label{lem:UniformConvergence}
For any $n\in \BN$, there exists an $r_n>0$ such that the product
\[\prodprm_{a\in A} (1+f_a(u))\]
converges uniformly for $(u,\omega')\in B(0,r_n)\times \Omega^{r-1}_n$.
\end{Lem}
\begin{Proof}
By definition $\varphi_a(X)=aX\prodprm_{\varphi_a(\alpha)=0}\left(1-\frac{X}{\alpha}\right)$, 
and so $1+f_a(u)=\prodprm_{\varphi(\alpha)=0}\left(1-\alpha u\right)$. In these products, $\alpha$ 
runs through the values $e(a^{-1}\omega' A^{r-1})$.
We claim that these values are bounded from above by some $d$ that is independent 
of $a$ and $\omega'\in\Omega^{r-1}_n$, but may depend on $n$. In the definition of 
$e(z)=z\prodprm \left(1-\frac{z}{\lambda}\right)$ we can split the product into three factors: 
those where $\abs{\lambda}<\abs{z}$, those where $\abs{\lambda}=\abs{z}$ and those 
where $\abs{\lambda}>\abs{z}$. The factors where $\abs{z}<\abs{\lambda}$ have 
absolute value 1 and the factors where $\abs{z}=\abs{\lambda}$ have absolute value 
less than or equal to 1. Thus
\begin{equation}\label{eq:BoundExponent}
 \abs{e(z)}\le \abs{z}\prodprm_{\abs{\lambda}<\abs{z}}\abs{\frac{z}{\lambda}}.
\end{equation}
Each $\lambda\in \omega' A^{r-1}$ is of the form $\lambda = a_1\omega_1+\cdots+a_{r-1}\omega_{r-1}$, (each $a_i\in A$), so
\[ \abs{\lambda}\ge \abs{\omega}_i \ge q^{-n}\abs{\omega}.\]
The first inequality follows since $a_1\omega_1+\cdots+a_{r-1}\omega_{r-1}$ is a 
$\Finf$-linear equation with largest coefficient greater or equal to 1, and the 
second follows from the definition of $\Omega^{r-1}_n$. It follows that each factor 
of \eqref{eq:BoundExponent} satisfies 
$\abs{\frac{\alpha}{\lambda}}\le \abs{\frac{\omega}{\lambda}}\le q^n$, 
since we may take $\abs{\alpha}\le \abs{\omega}$. It remains to bound the number 
of factors. Since $\lambda = a_1\omega_1+\cdots+a_{r-1}\omega_{r-1}$ and the 
smallest $\omega_i$ is at least as large as $q^{-n}$ times the largest $\omega_j$, 
the $a_i$ are bounded by $q^n$. Therefore the number of $\lambda$ is bounded by 
$q^{(r-1)n}$ which depends only on $n$. That proves the claim.\\

Now let $s_i$ ($i=1,q,\ldots,q^{(r-1)\deg a}$) be the symmetric polynomials in 
the set of $\alpha$ for which $\varphi_a(\alpha)=0$, so that
\[f_a(u) = \sum_{i=1}^{q^{(r-1)\deg a}}s_{q^{(r-1)\deg a}-i}u^{q^{(r-1)\deg}-i}.\]
Since, for each $\alpha$ we have $\abs{\alpha}<d$, depending only on $n$, we have
$\abs{s_j}\le d^j$.
For any $\epsilon<1$, suppose that $|u|<\epsilon/d$. Then 
$|s_ju^j|\le d^j\cdot (\epsilon/d)^j\le \epsilon^j$ for each $j=q^{(r-1)\deg a}-i$. It follows that
\[ \left|\sum_{i=1}^{q^{(r-1)\deg a}}s_{q^{(r-1)\deg a}-i}u^{q^{(r-1)\deg}-i}
  \right|\le \epsilon^{q^{(r-1)\deg a}-q^{(r-1)\deg a-1}}. \]
As $|a|\to \infty$, this tends to 0, so the product converges uniformly on $B(0,\epsilon/d)\times\Omega^{r-1}_n$.
\end{Proof}

\begin{Thm}\label{Thm:prodformula}
The $u$-expansion of the Drinfeld discriminant function $\Delta$ 
is given by the product
\[ \Delta(\omega) = -\Delta'(\omega')^q u^{q-1}\prodprm_{a\in A}(1+f_a(u))^{q^r-1}.\]
\end{Thm}

\begin{Rem}
If $r=2$, then we set $\Delta'=1$ to be consistent with the product formula determined by Gekeler. The cautious reader may note that the factor of $\period^{q^2-1}$ from Gekeler's formula is missing. This is a consequence of having normalized our elements $\omega\in\Omega^r$ so that $\omega_r=\period$ instead of 1.
\end{Rem}

\medskip

\begin{Cor}
If $A_+$ denotes the set of monic polynomials in $A$, then
\[ \Delta(\omega) = -\Delta'(\omega')^q u^{q-1}\prodprm_{a\in A_+}(1+f_a(u))^{(q^r-1)(q-1)}.\]
\end{Cor}
\begin{Proof}
Every polynomial can be uniquely written in the form $cp$ with $c\in \BF_q^\times$ and 
$p\in A_+$. All that remains is to note that by definition $f_{cp}(u)=f_{p}(u)$.
\end{Proof}

\begin{Proofof}{Theorem \ref{Thm:prodformula}}
The proof is very similar to Gekeler's proof in the rank 2 case \cite{GekelerProd}, so we give only the 
essential details. During this proof, let $\Lambda$ be the lattice $\omega_1A+\cdots+\omega_r A$, and 
$\Lambda'$ be the lattice $\omega_2A+\cdots+\omega_r A$. The following formulae from \cite{GekelerProd} 
still hold:

\begin{equation}
\Delta(\omega)=t\prodprm_{\alpha\in(t^{-1}A/A)^r} e_{\Lambda}(\omega\alpha)^{-1};
\label{eq:DeltaDef}
\end{equation}

\begin{equation}
e_{\Lambda}(X)=e_{\Lambda'}(X)\prodprm_{a\in A}\frac{e_{\Lambda'}(X)+e_{\Lambda'}(a\omega_1)}{e_{\Lambda'}(a\omega_1)};
\label{eq:ExponentialProduct}
\end{equation}
\begin{equation}
\Delta'(\omega')\cdot\!\!\!\!\!\!  \prod_{\substack{z\in\Cinf\\\varphi_t(z)=\varphi_t(z_0)}}(X-z)=\varphi_t(X-z_0),
\label{eq:DeltaDrinfeld}
\end{equation}
where $\varphi$ is the Drinfeld module associated to the lattice $\Lambda'$, and $z_0\in\Cinf$ is fixed. (They are equations (5), (6) and (4), respectively in \cite{GekelerProd}.)
From \eqref{eq:DeltaDef} we compute
\begin{align*}
\Delta(\omega)&=t\prodprm_{\alpha\in(t^{-1}A/A)^r} e_{\Lambda}(\omega\alpha)^{-1}\\
&\stackrel{\eqref{eq:ExponentialProduct}}{=\joinrel=} t\prodprm_{\alpha\in(t^{-1}A/A)^r}\left(e_{\Lambda'}(\omega\alpha)^{-1}\prodprm_{a\in A}\frac{e_{\Lambda'}
(a\omega_1)}{e_{\Lambda'}(a\omega_1)+e_{\Lambda'}(\omega\alpha)}\right)\\
&= t\left(\prodprm_{\alpha\in(t^{-1}A/A)^r}e_{\Lambda'}(\omega\alpha)^{-1}\right)\cdot 
\prodprm_{a\in A}\prodprm_{\alpha\in(t^{-1}A/A)^r}\frac{e_{\Lambda'}
(a\omega_1)}{e_{\Lambda'}(a\omega_1)+e_{\Lambda'}(\omega\alpha)}
\end{align*}
In the last product we note that adding the factor corresponding to $\alpha=0$, changes nothing, since it is equal to $\frac{e_{\Lambda'}(a\omega_1)}{e_{\Lambda'}(a\omega_1)}=1$. Thus
\[  \Delta(\omega) = t\left(\prodprm_{\alpha\in(t^{-1}A/A)^r}
e_{\Lambda'}(\omega\alpha)^{-1}\right)\cdot 
\prodprm_{a\in A}\prod_{\alpha\in(t^{-1}A/A)^r}\frac{e_{\Lambda'}
(a\omega_1)}{e_{\Lambda'}(a\omega_1)+e_{\Lambda'}(\omega\alpha)} \]
First we compute the internal product in the second factor. We break it up into two parts:

\begin{align}
\prod_{\alpha\in (t^{-1}A/A)^r}e_{\Lambda'}(a\omega_1)&=
e_{\Lambda'}(a\omega_1)^{q^r}\nonumber \\
&= \varphi_a(e_{\Lambda'}(\omega_1))^{q^r}  \nonumber\\ 
&\stackrel{\eqref{eq:UaInTermsOfHa}}{=\joinrel=} \Delta'_a(\omega')^{q^r}(1+f_a(u))^{q^r},
\label{eq:DiscNumerator}
\end{align}
and
\begin{align}
\prod_{\alpha\in(t^{-1}A/A)^r} \left(e_{\Lambda'}(a\omega_1)+e_{\Lambda'}(\alpha\omega)\right)^{-1}
&\stackrel{\eqref{eq:DeltaDrinfeld}}{=\joinrel=} \prod_{c\in (t^{-1}A/A)}\Delta'(\omega')\varphi_t\left(e_{\Lambda'}(a\omega_1)+e_{\Lambda'}(c\omega_1)\right)^{-1}\nonumber \\
&=\Delta'(\omega')^q\prod_{c\in\BF_q}\varphi_{at+c}(e_{\Lambda'}(\omega_1))^{-1}\nonumber \\
&\stackrel{\eqref{eq:UaInTermsOfHa}}{=\joinrel=}\Delta'(\omega')^q\prod_{c\in\BF_q}\Delta'_{at+c}(\omega')(1+f_{at+c}(u)).
\label{eq:DiscDenominator}
\end{align}
Together \eqref{eq:DiscNumerator} and \eqref{eq:DiscDenominator} give
\[
\prodprm_{a\in A}\frac{e_{\Lambda'}(a\omega_1)}{e_{\Lambda'}(a\omega_1)+e_{\Lambda'}(\omega\alpha)}
 = \prodprm_{a\in A} \frac{\Delta'_a(\omega')^{q^r}(1+f_a(u))^{q^r}}{\Delta'(\omega')^q\prod_{c\in\BF_q}\Delta'_{at+c}(\omega')(1+f_{at+c}(u))}\nonumber 
\]
If $\deg a\ge 1$ we use \eqref{eq:DeltaA} to obtain 
$\Delta'_{at+c}(\omega')=\Delta'_a(\omega')^{q^{r-1}}\Delta'(\omega')$ after which 
the previous equation simplifies to 
\begin{equation}
\prodprm_{a\in A}\frac{e_{\Lambda'}(a\omega_1)}{e_{\Lambda'}(a\omega_1)+e_{\Lambda'}(\omega\alpha)} 
=\prodprm_{a\in A} \frac{(1+f_a(u))^{q^r}}{\prod_{c\in\BF_q}(1+f_{at+c}(u))}, 
\label{eq:DiscNumeratorAndDenominator}
\end{equation}
If $\deg a=0$, it becomes instead
\begin{equation}
\prodprm_{a\in\BF_q}\frac{1}{\prod_{c\in\BF_q}(1+f_{at+c}(u))}.
\label{eq:DiscProdDegree0}
\end{equation}

We now claim that
\[\prod_{\substack{a\in A\\ \deg a\ge 1}}(1+f_a(u))=\prodprm_{a\in A}\prod_{c\in\BF_q}(1+f_{at+c}(u)).\]

Indeed, these products are of the same factors over the same index set, so the result follows by 
uniform convergence of the products involved (Lemma \ref{lem:UniformConvergence}). Since 
$f_a(u)=0$ when $\deg a=0$, 
\eqref{eq:DiscNumeratorAndDenominator} simplifies to
\begin{equation}
\prodprm_{a\in A} (1+f_a(u))^{q^r-1}.
\label{eq:DiscNumeratorAndDenominatorSimplified}
\end{equation}

Lastly, we compute 
\[ \prodprm_{\alpha\in (t^{-1}A/A)^r}e_{\Lambda'}(\omega\alpha),\] 
 again by breaking it into two parts. Let $\alpha=(\alpha_1,\alpha')$ where 
$\alpha_1\in (t^{-1}A/A)$ and $\alpha'\in (t^{-1}A/A)^{r-1}$ and distinguish
the two cases where $\alpha_1=0$ and where $\alpha_1\ne 0$. Using 
\eqref{eq:DeltaDrinfeld}, we get
\begin{align}
\prodprm_{\alpha_1\in(t^{-1}A/A)}\prod_{\alpha'\in (t^{-1}A/A)^r}e_{\Lambda'}(\omega\alpha)&=
\prodprm_{\alpha_1\in(t^{-1}A/A)}\Delta'(\omega')\varphi_t(e_{\Lambda'}(\omega_1\alpha_1))^{-1}\nonumber\\
&= \Delta'(\omega')^{q-1} \prodprm_{c\in\BF_q}cu\nonumber\\
&=-\Delta'(\omega')^{q-1}u^{q-1}.
\label{eq:FirstFactorPart1}
\end{align}
The factor with $\alpha_1=0$ is
\begin{equation}
\prodprm_{\alpha\in (t^{-1}A/A)^{r-1}}e_{\Lambda'}(\omega'\alpha')=\Delta'(\omega')
\label{eq:FirstFactorPart2}
\end{equation}

Putting together \eqref{eq:DiscNumeratorAndDenominatorSimplified}, \eqref{eq:FirstFactorPart1} 
and \eqref{eq:FirstFactorPart2} we obtain the Theorem.
\end{Proofof}

\begin{Rem}
In this note we have restricted to the simplest case where $A=\BF_q[t]$ and $\Delta$ 
is the discriminant function that arises from the polynomial $\varphi_t(X)$. In 
\cite[\S VII.4]{GekelerDMC}, Gekeler extended this result to general Drinfeld 
rings $A$ where the discriminant function is the one associated to $\varphi_a(X)$ 
with $a\in A$. The theory of higher rank Drinfeld modular forms as developed in 
\cite{BBP} allows for this generality and it is probably not hard to extend this 
result to the more general setting. 

In \cite{Hamahata2} Hamahata has already generalized his form 
of the product expansion to general Drinfeld rings $A$. 
\end{Rem}

\noindent\textbf{Acknowledgements.} The author wishes to thank Florian Breuer for 
many discussions during the author's PhD and many comments and suggestions in 
preparation of this article and Richard Pink who initiated 
the project on higher rank Drinfeld modular forms  with Florian Breuer. Without them this work
would not have been possible.


\end{document}